\documentclass[oneside,a4paper,12pt]{article}

\usepackage{geometry} 
\usepackage{amsmath}
\usepackage{amsthm}
\usepackage{amsfonts,amssymb}
\usepackage{fancyhdr} 
\usepackage{mathrsfs}  
\usepackage{graphicx}
\usepackage[all]{xy}

\usepackage{times}
\usepackage{cases}
\usepackage{color}
\usepackage[colorlinks,linkcolor=blue,anchorcolor=blue,citecolor=blue]{hyperref}

\usepackage{float} 
\usepackage{subfig}
\usepackage{amsthm}

\newtheorem{theorem}{Theorem}[section]
\newtheorem{lemma}[theorem]{Lemma}
\newtheorem{proposition}[theorem]{Proposition}
\newtheorem{definition}[theorem]{Definition}
\newtheorem{corollary}[theorem]{Corollary}

\newtheorem{claim}[theorem]{Claim}

\theoremstyle{remark}
\newtheorem{remark}[theorem]{Remark}

\theoremstyle{remark}

\numberwithin{equation}{section}

\begin{document}

\title{\textbf{Angle structure on general hyperbolic $3$-manifolds}}

\author{\medskip Huabin Ge, Longsong Jia, Faze Zhang}

\date{}

\maketitle
\begin{abstract}

 Let $M$ be a non-compact hyperbolic $3$-manifold with finite volume and totally geodesic boundary components. By subdividing mixed ideal polyhedral decompositions of $M$, under some certain topological conditions, we prove that $M$ has an ideal triangulation which admits an angle structure.
\end{abstract}

\maketitle




\section{Introduction}

Thurston's geometrization conjecture (proved by Perelman~\cite{perelman1}-\cite{perelman3}, or see \cite{CaoZhu},\cite{KleiLott},\cite{MorTian} for more details) tells us that any $3$-manifold can be decomposed via J-S-J splitting into $8$ distinct geometric pieces, the majority of which are hyperbolic. A $3$-manifold with a hyperbolic geometric structure is termed as a hyperbolic $3$-manifold. Thus, determining the existence of a hyperbolic structure on a given $3$-manifold is a significant research topic. Thurston~\cite{Thurston} introduced the concept of ideal triangulations on $3$-manifolds and proved that a solution of the hyperbolic gluing equations on the ideal triangulation can yield a hyperbolic structure on the $3$-manifold. However, solving the hyperbolic gluing equations on general ideal triangulations is extremely challenging. In the 1990s, Casson \cite{Lackenby-1} and Rivin \cite{Rivin} discovered a powerful technique for solving Thurston's gluing equations. By introducing the concept of angle structures on ideal triangulations, they proved that if a maximal volume angle structure exists, it provides solutions to Thurston's gluing equations, thereby yielding a hyperbolic structure on the $3$-manifold. The converse is also true, and one can find a self-contained exposition of Casson-Rivin's program in \cite{FuGue}.

An \emph{angle structure} on an ideal triangulation refers to assigning real numbers (called dihedral angle) in the interval $(0,\pi)$ at each edge of each ideal tetrahedron such that the sum of the dihedral angles at each ideal vertex is $\pi$ and the sum of dihedral angles around each edge is $2\pi$. Kang-Rubinstein~\cite{KR} extended the real numbers assigned at each edge into $[0,\pi]$, and called this structure as \emph{semi-angle structure}. If the dihedral angles are only taken values $0$ or $\pi$, then the structure is termed as \emph{taut structure} by Lackanby~\cite{Lac}. Luo-Tillmann~\cite{LT} further extended the assignment of dihedral angles to the whole real numbers $\mathbb{R}$, which is defined as \emph{generalized angle structure}.

Conversely, to study the geometry and topology of hyperbolic $3$-manifolds, a natural approach is to investigate their geometric ideal tetrahedra decompositions. The existence of geometric ideal triangulations are very important for studying hyperbolic 3-manifolds. On the one hand, Thurston's proof of the famous Hyperbolic Dehn Filling Theorem takes it for granted that each cusped hyperbolic 3-manifold has a geometric ideal triangulation, see \cite{BenPet}, \cite{Martelli}, \cite{Petronio} for example. On the other hand, geometric ideal triangulations ensures that all three-dimensional hyperbolic manifolds can be constructed from ideal tetrahedra with the help of Thurston's gluing equations. The first breakthrough on geometric ideal decomposition is due to Epstein-Penner~\cite{EP}. They introduced canonical geometric ideal polyhedral decompositions on cusped hyperbolic manifolds with finite volume. Similarly, Kojima~\cite{Kojima} obtained truncated hyperideal polyhedral decompositions on compact hyperbolic $3$-manifolds with totally geodesic boundary. Recently, we \cite{GJZ} obtained mixed ideal polyhedral decompositions on non-compact hyperbolic $3$-manifolds with finite volume and totally geodesic boundary. 
In all known examples, a cusped hyperbolic manifold $M$ admits at least one geometric ideal triangulation, but whether this holds true for arbitrary cusped hyperbolic $3$-manifold $M$ is still an open question. As was pointed by Gu\'{e}ritaud and Schleimer in \cite{Gueritaud}, ``it is a difficult problem in general. General results are known only when $M$ is restricted to belong to certain classes of manifolds: punctured-torus bundles,
two-bridge link complements, certain arborescent link complements and related objects, or covers of any of these spaces" (see \cite{Aki}-\cite{Aki-SWY-2}, \cite{Gue-1}-\cite{Gue-2}, \cite{Ham-P}, \cite{Jorgen}, \cite{Lackenby}, \cite{Nimer} for instance).
For compact three-dimensional manifolds with boundary, Feng-Ge-Hua \cite{Feng-Ge-Hua} established a deep connection between the combinatorial Ricci flow and geometric ideal triangulations, and proved that any topological ideal triangulation with degree greater than or equal to 10 are geometric. Although it was difficult to obtain geometric triangulations directly, Luo-Schleimer-Tillmann \cite{Luo-S-T} proved the virtual existence of geometric triangulations. Furthermore, Futer-Hamilton-Hoffman \cite{Futer-H-H} proved that there are infinitely many virtual geometric triangulations.

To a certain extent, the existence of angle structure is weaker than that of geometric triangulations. Hence, a quite natural question arises: for any hyperbolic $3$-manifold $M$, are there topological triangulations which admits an angle structure? Luo-Tillmann~\cite{LT} gave a sufficient and necessary condition for the existence of angle structure on cusped hyperbolic $3$-manifolds using the theory of normal surface. By using Epstein-Penner's  ideal geometric polyhedral decompositions and the results of Luo-Tillmann~\cite{LT}, Hodgson-Rubinstein-Segerman~\cite{HRS} proved that there exists an ideal triangulation which admits an angle structure on cusped hyperbolic $3$-manifolds with a certain topological condition. By using Kojima's hyperideal (truncated) polyhedral decompositions, Qiu-Zhang-Yang~\cite{QZY} obtained a hyperideal (truncated) triangulation which admits an angle structure on compact hyperbolic $3$-manifolds with totally geodesic boundary. In this paper we consider a hyperbolic $3$-manifold $M$ with both cusps and totally geodesic boundaries. By Moise's work \cite{Moise}, there is always an ideal triangulation $\mathcal T$ on $M$. We want to know whether $\mathcal T$ admits an angle structure.
Based on our previous work on mixed ideal polyhedral decompositions \cite{GJZ} and Luo-Tillmann's arguments in \cite{LT}, we get a natural ideal triangulation $\mathcal T$ (see Corollary~\ref{cor}) on $M$ directly without using Moise's work \cite{Moise}. Moreover, under some topological conditions, $\mathcal T$ admits an angle structure. To be precise, we have:

\begin{theorem} \label{main2}
Let $M$ be a non-compact finite-volume hyperbolic $3$-manifold with totally geodesic boundary. Denote $\overline{M}$ by the compact $3$-manifold with boundary, with each torus (or Klein bottle) boundary component corresponds to a cusp of $M$. In other words, after subtracting the torus (or Klein bottle) boundary components, $\overline{M}$ is homeomorphic to $M$. If $H_{1}(\overline{M};\mathbb Z_{2})\rightarrow H_{1}(\overline{M},\partial \overline{M};\mathbb Z_{2})$ is the zero map, then there is a topological ideal triangulation on $M$ such that it admits an angle structure.
\end{theorem}

To prove Theorem \ref{main2}, the polyhedral decomposition theory \cite{GJZ}, which parallels \cite{EP} and \cite{Kojima}, will play a crucial role. In addition, we follow the spirits in Luo-Tillmann \cite{LT}, Kang-Rubinstein~\cite{KR} and Hodgson-Rubinstein-Segerman~\cite{HRS}. We outline the proof in three steps. Step 1, using our geometric polyhedral decompositions \cite{GJZ}, we first decompose $M$ into hyperbolic ideal or truncated polyhedra. Then we triangulate each polyhedron. The triangulation of the common faces of adjacent polyhedra may not necessarily match. On the mismatched faces, by inserting flat ideal tetrahedra, we are able to achieve an ideal triangulation $\mathcal T$. Step 2, following Kang-Rubinstein~\cite{KR} and Luo-Tillmann~\cite{LT}, there is a strong relationship between the angle structure of $(M,\mathcal T)$ and the theory of normal surfaces in $\mathcal T$. On each individual $\sigma_{i}$ of $\mathcal T$, there exist four types of normal triangles and three types of normal quadrilaterals, each of which is characterized by coordinates. These coordinates satisfy a system of \emph{compatibility equations} with solution space denoted by $\mathcal C(M,\mathcal T)$. By introducing the \emph{generalized Euler characteristic function} $\chi^{*}$ defined on $\mathcal C(M,\mathcal T)$ and using Farkas's lemma, we derive a sufficient combinatorial condition for the existence of angle structures, i.e.

\begin{proposition} \label{main1}
Let $M$ be a non-compact hyperbolic $3$-manifold with finite volume and totally geodesic boundary. There exists an ideal triangulation $\mathcal T$ on $M$, and if $\chi^{*}(s)<0$ for all $s\in \mathcal C(M, \mathcal T)$ with all quadrilateral coordinates non-negative and at least one quadrilateral coordinates positive, then $(M,\mathcal T)$ admits an angle structure.
\end{proposition}
\noindent Step 3, by excluding cases that do not satisfy the combinatorial conditions in Proposition \ref{main1}, we get the topological conditions in Theorem \ref{main2} for the existence of angle structures.

~

The paper is organized as follows. We give some basic notions in Section \ref{2}, including the polyhedral decomposition theory, ideal triangulation, angle structure, normal surface and Farkas's lemma. In Section \ref{3}, we derive an ideal triangulation directly and prove Proposition \ref{main1}. In Section \ref{4}, we prove Theorem \ref{main2}.


~

\noindent
\textbf{Acknowledgements:}
The authors are very grateful to Professor Ruifeng Qiu, Feng Luo, Tian Yang for many discussions on related problems in this paper. The first two authors would like to thank Professor Gang Tian for his long-term support and encouragement. Huabin Ge is supported by NSFC, no.12341102 and no.12122119. Faze Zhang is
supported by NSFC, no.12471065.


\section{Preliminaries}\label{2}

\subsection{Polyhedral decomposition theory}
Let us recall the definitions of the ideal tetrahedra and partially truncated tetrahedra first.

\begin{figure}[htbp]
\centering
\includegraphics[scale=0.27]{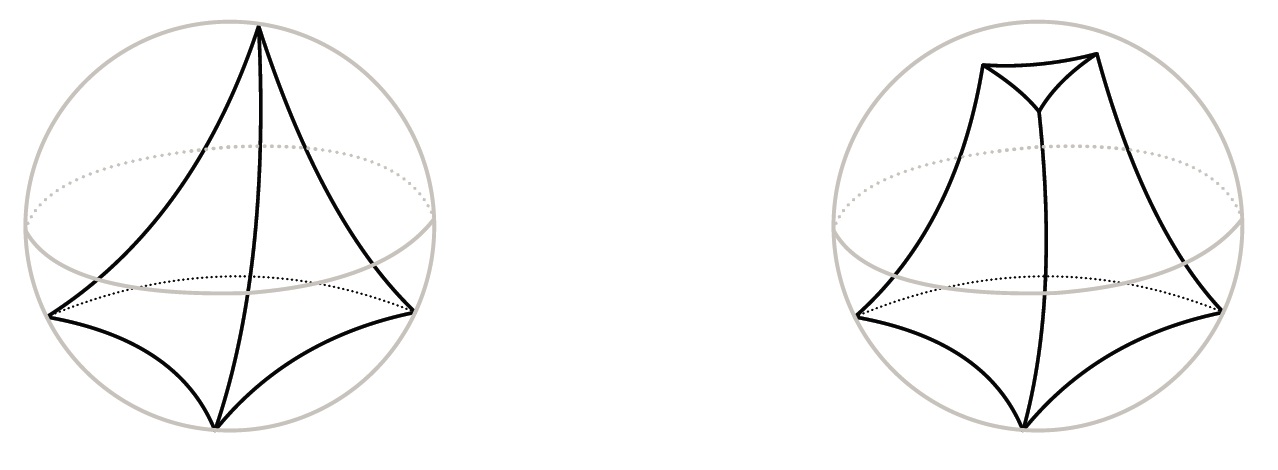}
\caption{an ideal tetrahedron (left) / a type 1-3 partially truncated tetrahedron (right)}
\label{Fig1}
\end{figure}

\begin{definition}
An \textbf{ideal tetrahedron} in $\mathbb {H}^{3}$ is a finite-volume region in $\mathbb {H}^{3}$ bounded by four geodesic faces. Any two faces intersect each other,
and any three faces intersect at the infinite boundary of $\mathbb {H}^{3}$. See the left picture of Figure \ref{Fig1} for example.
\end{definition}

Denote $\mathbb{H}^3\subset RP^{3}$ by the open unite ball representing the hyperbolic 3-space  via the projective Klein model. Following the terminology in \cite{BB}, we have
\begin{definition}
A \textbf{hyperideal tetrahedron} $\hat{P}$ is a non-compact tetrahedron in $\mathbb{H}^3$ which is just the intersection of $\mathbb{H}^3$ with a projective tetrahedron $\tilde{P}$ whose vertices are all outside $\mathbb{H}^3$ and whose edges all meet $\mathbb{H}^3$.
\end{definition}

Given that the hyperideal tetrahedron has been extensively studied and used by several mathematicians, we would like to emphasize that the definition of a hyperideal tetrahedron varies from author to author. Some authors require that all four vertices of a hyperideal tetrahedron must be hyperideal vertices, i.e. points outside the closure of $\mathbb{H}^3$, e.g. \cite{Luo-Y}. It is called \emph{strictly hyperideal} by Schlenker \cite{Schlenker}.

\begin{definition}\label{def1}
A \textbf{hyperideal tetrahedron of 1-3 type} $\hat{P}$ means that there are only one vertex $v$ of $\tilde{P}$ lies outside of the closure of $\mathbb{H}^3$ and the other three vertices are all lie in $\partial \mathbb {H}^{3}$. The truncation of a type 1-3 hyperideal tetrahedron $\hat{P}$ at the vertex $v$ is defined as cutting off the thick end towards $v$ from $\hat{P}$ by a geodesic plane which is perpendicular to the three adjacent geodesic faces at $v$. After the truncation of $\hat{P}$ at  $v$, the resulting polyhedron $P$ is called a \textbf{partially truncated tetrahedron of 1-3 type}. The unique face of $P$ induced by the truncation is called an external face and the other faces of $P$ are called internal faces. See the right picture of Figure \ref{Fig1} for example.
\end{definition}

The concept of ideal tetrahedron, hyperideal tetrahedron and partially truncated tetrahedra of 1-3 type can be generalized to polyhedra case easily.

\begin{definition}
An \textbf{ideal polyhedron} in $\mathbb {H}^{3}$ is a finite-volume polyhedral region bounded by several geodesic faces. Any two faces intersect each other, and any three faces intersect at the infinite boundary of $\mathbb {H}^{3}$.
\end{definition}

\begin{definition}
A \textbf{hyperideal polyhedron} $\hat{P}$ is a non-compact polyhedron in $\mathbb{H}^3$ which is just the intersection of $\mathbb{H}^3$ with a projective polyhedron $\tilde{P}$ whose vertices are all outside $\mathbb{H}^3$ and whose edges all meet $\mathbb{H}^3$.
\end{definition}

\begin{definition}
A \textbf{hyperideal polyhedron of 1-k type} $\hat{P}$ is defined as there is exactly one vertex $v$ (called hyperideal vertex) of the corresponding $\tilde{P}$ lies outside of $\mathbb{H}^3\cup \partial \mathbb {H}^{3}$ and the other $k$ vertices (called ideal vertices) are all lie in $\partial \mathbb {H}^{3}$. After the truncation of $\hat{P}$ at the hyperideal vertex $v$, the resulting polyhedron $P$ is defined as a \textbf{partially truncated polyhedron of 1-k type}. The unique face of $P$ induced by the truncation is called an external face and the other faces of $P$ are called internal faces.
\end{definition}

Similarly, one can define the truncation of a hyperideal tetrahedron (or polyhedron resp.), which is truncated at all hyperideal vertices, i.e. points outside the closure of $\mathbb{H}^3$. After the truncating, one get a partially truncated tetrahedron (or polyhedron resp.).

Similar to the role of Epstein-Penner's decomposition \cite{EP} played in \cite{HRS}, we developed a theory of decomposition of general hyperbolic 3-manifolds, which will paly essential role in the proof of our main theorem. The theory says any hyperbolic 3-manifold admits a geometric polyhedral decomposition.
\begin{theorem}[Ge-Jiang-Zhang, \cite{GJZ}]
\label{the}
Let $M$ be a volume finite, non-compact, complete hyperbolic $3$-dimensional manifold with totally
geodesic boundary. Then $M$ admits a mixed geometric ideal polyhedral decomposition $\mathcal{D}$ such that each cell is either an ideal polyhedron or a partially truncated polyhedron with only one hyperideal vertex.
\end{theorem}

Theorem \ref{the} may be considered as a combination of Epstein-Penner's decomposition \cite{EP} and Kojima's decomposition \cite{Kojima}. However, to some extent, our conclusion is slightly stronger, and only ideal polyhedra and 1-$k$ type polyhedra appear in our decomposition. It is worth emphasizing that there are several advantages to choosing type 1-3 tetrahedra and type 1-$k$ polyhedra. Firstly, for polyhedral decomposition, it is more suitable for the algorithm. Secondly, it is very suitable for the proof of the main theorem in this paper, which can highlight the essence of the proof, and can also greatly simplify the expression, making the details clearer and more readable. See Section \ref{3}-\ref{4} for details.

\subsection{Ideal triangulation}


Strictly speaking, the ideal tetrahedron, as well as the hyperideal tetrahedron, have no real existent vertices. But for the sake of convenience, we often think that they all have four vertices. For example, a truncated type 1-3 tetrahedral $\sigma$ has three real vertices, which are the three vertices of a hyperbolic triangle obtained by truncating a strictly hyperideal vertex, and four non-existent vertices, of which three are ideal vertices and one is a strictly hyperideal vertex. But we often think that these four unreal vertices are all the vertices of $\sigma$, while those three real vertices are not considered vertices of $\sigma$.



Additionally, it should be noted that the geometry of ideal tetrahedra, hyperideal tetrahedra, and truncated tetrahedra is hyperbolic. For each ideal tetrahedron, it corresponds to a topological tetrahedron with no geometry (i.e., on the boundary of a topological 3-ball, four marked vertices are selected, and each of the two marked points is connected to a marked edge, which does not intersect each other). It can even be considered equivalently that these topological tetrahedra are combinatorial tetrahedra with only combinatorial structures, where combinatorial structures refer to vertices, edges, faces, and connection relationships between them. We can also think of it as corresponding to a type of Euclidean tetrahedra, which may not be the same in shape, but have the same combinatorial structure. Similarly, each truncated hyperbolic tetrahedron corresponds to a combinatorial tetrahedron and a type of truncated Euclidean tetrahedron. 
For example, the following Figure \ref{Fig-correspond} depicts the correspondence between a tetrahedron with only one truncated vertex (left) and a tetrahedron with four truncated vertices (right). If a (truncated) Euclidean tetrahedron $\sigma$ corresponds to a (truncated) hyperbolic tetrahedron according to the above rules, then a vertex $v$ in $\sigma$ is called an \emph{ideal vertex} or a \emph{hyperideal vertex}, defined according to its corresponding vertex type in the (truncated) hyperbolic tetrahedron.

\begin{figure}[htbp]
\centering
\includegraphics[scale=0.18]{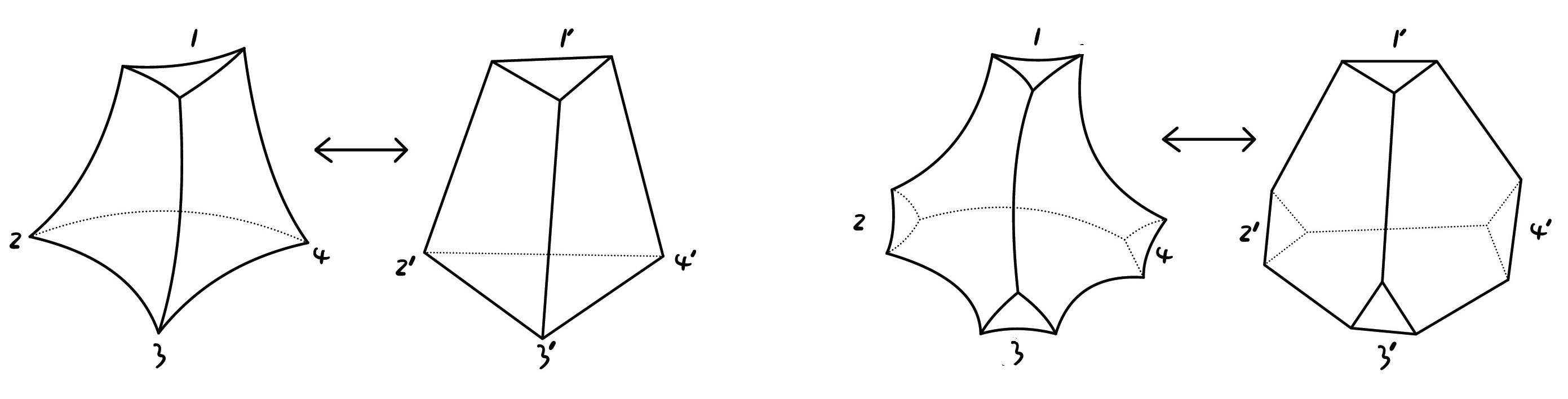}
\caption{correspondence of tetrahedra}
\label{Fig-correspond}
\end{figure}

The surface generated by truncating a vertex is called an \emph{external surface}. The remaining faces are called \emph{internal faces}, or in other words, faces containing at least one ideal vertex are called internal faces. In Euclidean tetrahedra and combinatorial tetrahedra, external and internal faces can be defined similarly. The edges in the external faces are called the \emph{external edges}, and the remaining edges are called the \emph{internal edges}. Equivalently, the edges that contain at least one ideal vertex are called the internal edges. Thus, each tetrahedron, whether it is an ideal tetrahedron or a truncated tetrahedron, always has exactly six internal edges. Below we provide the precise meaning of an ideal triangulation.

\begin{definition}
\label{Def-ideal-triangulation}
Given a topological three dimensional manifold $M$ with both topological cusps and boundaries whose connected components all have negative Euler characteristic numbers. For example, $M$ is derived from a volume finite, non-compact, complete hyperbolic $3$-manifold with totally geodesic boundary, after forgetting the hyperbolic structure.
Let $X=\sqcup_{i=1}^n\sigma_{i}$ be a topological space with each $\sigma_{i}$ either an Euclidean tetrahedron or a truncated Euclidean tetrahedron. Let $\Phi$ be a collection of affine homeomorphisms between the internal faces, and $X^{(0)}$ be the set of ideal vertices in $X$. If $(X\backslash X^{(0)})/\Phi$ is homomorphic to $M$, then the cell decomposition $\mathcal{T}=X\backslash X^{(0)}$ is called an \textbf{ideal triangulation} of $M$.
\end{definition}

Note that the quotient of external faces in $X$ constitutes a triangulation of $\partial M$. The $1$-simplices ($2$-simplices resp.) in $M\setminus\partial M$ are called internal edges (faces resp.) of $\mathcal{T}$.

\subsection{Angle structure}
\label{subsection:2.1}

By definition, an angle structure on a single topological, or combinatorial ideal tetrahedron $\sigma$ assigns to each edge of $\sigma$ a positive number, called \emph{dihedral angle}, such that the sum of the angles associated to the three edges meeting in a vertex is $\pi$ for each vertex of the tetrahedron. After a simple calculation, it can be found that the dihedral angles of opposite sides are equal. For a truncated tetrahedron with at least one hyperideal vertex, the definition of an angle structure is a little different. The main difference is that the sum of the dihedral angles associated to the three edges meeting in a strictly hyperideal vertex (or equivalently, a truncated face) is less than $\pi$. Combining the above two definitions, the concept of an angle structure can be extended to the most general 3-dimensional pseudo-manifolds with ideal triangulations, see \cite{Garoufa-H-R-S}, \cite{KR}, \cite{Lackenby-1}, \cite{Luo-3-flow}, \cite{LT}, \cite{Luo-Y} for instance.

\begin{definition}
\label{angle-structure}
Let $\mathcal T$ be an ideal triangulation on a topological 3-manifold $M$ with 3-simplices $\sigma_1, ... ,\sigma_n$. An \textbf{angle structure} on $(M,\mathcal T)$ is a function $\alpha$ that, for each internal edge $e_{ij}$ of $\sigma_i$ ($1\leq i\leq n$, $1\leq j \leq 6$), assigns a values $\alpha(e_{ij})\in(0,\pi)$ which is called the dihedral angle, so that:
\begin{enumerate}
\item for any internal edge $e$ in $\mathcal T$, the sum of all dihedral angles around $e$ is $2\pi$;
\item if $v$ is a vertex of some $\sigma_i$ with $1\leq i\leq n$, and is not in the external faces of $\sigma_{i}$, the sum of all dihedral angles adjacent to $v$ in $\sigma_{i}$ is $\pi$;
\item if $f$ is an external face of some truncated tetrahedron $\sigma_j$ with $1\leq j\leq n$, the sum of all dihedral angles adjacent to $f$ in $\sigma_j$ is less than $\pi$.
\end{enumerate}
\end{definition}

In the above definition, if all the dihedral angles are allowed to be taken in the closed interval $[0,\pi]$, it is said to be a \textbf{\emph{semi-angle structure}}. By the work of Bao-Bonahon \cite{BB} ( or see \cite{Luo-Y} and \cite{Schlenker}), any angle structure on $(M,\mathcal T)$ endows each single 3-simplex $\sigma_i$ with a hyperbolic geometry, making it an ideal tetrahedron or hyperideal truncated tetrahedron according to its type. However, it must be pointed out that due to the possible differences in edge lengths (even after choosing decorations), these hyperbolic tetrahedra generally cannot be glued together to obtain a hyperbolic structure on $M$.

\subsection{Normal surface and generalised Euler characteristic}\label{subsection:2.2}
Let $\mathcal T$ be an ideal triangulation on a topological 3-manifold $M$ with 3-simplex $\sigma_1, ... ,\sigma_n$.
\begin{definition}
Let $f$ be an internal face. A normal arc is an embedding of an arc $l$ into $f$ such that $\partial l$ lie on different internal edges of $f$. A normal disk in $\sigma_{i}$ ($1\leq i\leq n$) is an embedding of a disk $d$ into $\sigma_{i}$ ($1\leq i\leq n$) such that $d \cap \partial \sigma_{i}$ is either empty or consists of disjoint normal arcs.
\end{definition}

\begin{figure}[htbp]
\centering
\includegraphics[scale=0.35]{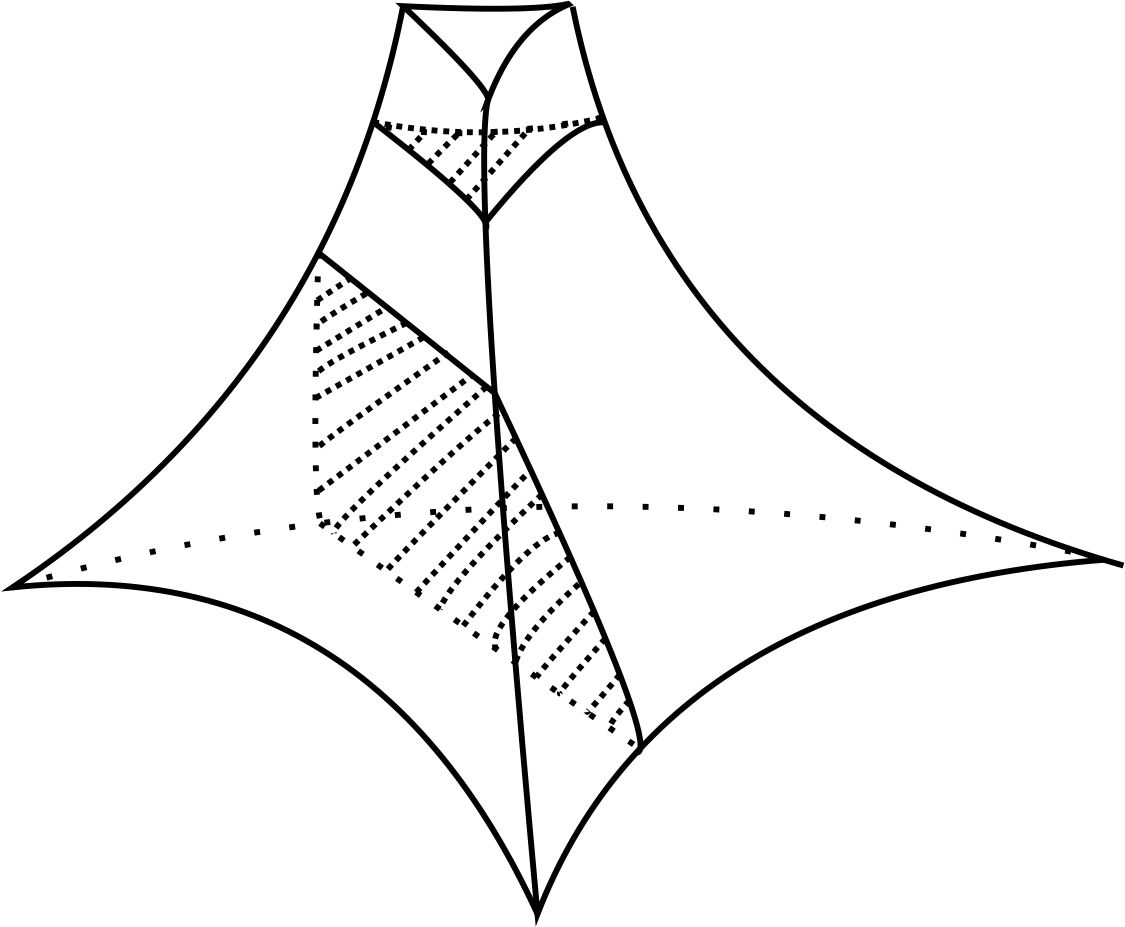}
\caption{normal disks}
\label{Fig4}
\end{figure}

For each $1\leq i\leq n$, there are seven types of normal disks on $\sigma_{i}$, of which 4 are normal triangles, corresponding to 4 vertices; 3 normal quadrilaterals, corresponding to 3 sets of opposing edges.
For example, Figure \ref{Fig4} shows a special normal triangle, as well as an example of a special normal quadrilateral.


\begin{definition}
A compact surface $F$ in $M$ is called \textbf{normal} with respect to $\mathcal T$, if $F \cap \sigma_{i}$ is either empty or consists of disjoint normal disks of $\sigma_{i}$ for each $1\leq i\leq n$.
\end{definition}

Recall $n$ is the number of all types of tetrahedra in the triangulation $\mathcal T$. We fix an ordering of all normal disc types $(q_1,..., q_{3n}, t_1,..., t_{4n})$ in $\mathcal T$, where $q_i$ denotes a normal quadrilateral type and $t_j$ a normal triangle type. For any given normal surface $F$, its normal coordinate $\overline{F}=(x_1,...,x_{3n},y_{1},...,y_{4n})$ is a vector in $\mathbb{R}^{7n}$, where $x_i$ is the number of normal discs of type $q_i$ in $F$, and $y_j$ is the number of normal discs of type $t_j$ in $F$. A properly embedded normal surface $F$ is uniquely determined up to normal isotopy by its normal coordinate (\cite{LT}, \cite{Tillmann}). $(x_1,...,x_{3n})$ is called the \emph{quadrilateral coordinate} of $F$.

For a normal disk $D$, let $b(D)$ be the number of normal arcs in $D\cap\partial M$. If $t$ is a normal triangle in $F$ meeting a particular tetrahedra at edges $e_i$ with edge valence $d_i$, $1\leq i \leq 3$, then its contribution to the Euler characteristic of $F$ is taken to be
$$\chi^{*}(t)=-\frac{1}{2}\big(1+b(t)\big)+\sum_{i=1}^3\frac{1}{d_i}.$$
If $q$ is a normal quadrilateral in $F$ meeting a particular tetrahedra at edges $e_i$ with edge valence $d_i$, $1\leq i \leq 4$, then its contribution to the Euler characteristic of $F$ is taken to be
$$\chi^{*}(q)=-\frac{1}{2}\big(2+b(q)\big)+\sum_{i=1}^4\frac{1}{d_i}.$$
Then the \emph{generalized Euler characteristic} function $\chi^{*}: \mathbb{R}^{7n}\rightarrow \mathbb{R}$ is defined as
\begin{equation}
\chi^{*}(x_1,...,x_{3n},y_{1},...,y_{4n} )=\sum_{i=1}^{3n}x_{i}\chi^{*}(q_{i})+\sum_{j=1}^{4n}y_{j}\chi^{*}(t_{j}),
\end{equation}
where $q_{i}$ and $t_{j}$ are the normal quadrilaterals and normal triangles in $F$ respectively. 
The generalized Euler characteristic function is linear, and coincides with the classical Euler characteristic for any embedded or immersed normal surface $F$ in $(M, \mathcal T)$. See Section 2.6, \cite{LT} for details.

There are some linear constraints between the normal coordinate of a normal surface $F$. Let $f$ be an internal face shared by two tetrahedra, say $\sigma_{i}$ and $\sigma_{j}$. Let $a$ be a normal isotopy class of a normal arc in $f$. There are two types of normal disks in $\sigma_{i}$, one is type $x_{i}$ for normal triangles, the other is type $y_{j}$ for normal quadrilaterals, such that when restricted to $f$, the corresponding normal arcs are of the same type as $a$. Similarly, there are two types of normal disks in $\sigma_{j}$, one is type $x_{i'}$ for normal triangles, the other is type $y_{j'}$ for normal quadrilaterals with similar properties. Thus, the normal coordinate $\overline{F}$ satisfies a linear equation at the face $f$ as follows, which is called the \emph{compatibility equation}:
\begin{equation}\label{equation1}
x_{i}+y_{j}=x_{i'}+y_{j'}.
\end{equation}
Denote $\mathcal C(M,\mathcal T)$ by the set of all vectors $(x_1,...,x_{3n},y_1,...,y_{4n})\in\mathbb{R}^{7n}$ satisfying (\ref{equation1}) at any internal face. For this linear space, Kang-Rubinstein~\cite{KR} once introduced a basis
$$\big\{W_{\sigma_{i}},\,W_{e_{j}}\,|\,i=1,...,n,\,j=1,...,m\big\},$$
where $n$ is the number of $3$-simplex in $\mathcal T$ and $m$ is the number of internal edges in $\mathcal T$. For concrete meanings of $W_{\sigma_{i}}$ and $W_{e_{j}}$, see \cite{KR} and \cite{LT}. Then for any $s\in \mathcal C(M,\mathcal T)$, there exists a unique coordinate $(\omega_{1},...,\omega_{n},z_{1},...,z_{m} )$ such that
\begin{equation}\label{equation2}
s=\sum_{i=1}^{n}\omega_{i}W_{\sigma_{i}}+\sum_{j=1}^{m}z_{j}W_{e_{j}}.
\end{equation}

Given a normal disk $d$ in some $\sigma_{i}$, if $d$ is a $k$-sided polygon, and intersects the internal edges $e_{i1},...,e_{ik}$ of some particular $ \sigma_{i}$, then its \emph{combinatorial area} $A(d)$ is defined to be
\begin{equation}
\sum_{j=1}^{k}\alpha(e_{ij})-(k-2)\pi.
\end{equation}

\begin{definition}\label{s-chi(A)}
If $s=\sum\limits_{i=1}^{n}\omega_{i}W_{\sigma_{i}}+\sum\limits_{j=1}^{m}z_{j}W_{e_{j}}$, then define
\begin{equation}\label{equation3}
\chi^{(A)}(s)=\frac{1}{2\pi} \sum_{t\in\bigtriangleup}x_{t}(s)A(t),
\end{equation}
where $\bigtriangleup$ is the set of all the normal triangles of $s$, and $x_{t}$ is the corresponding normal triangle coordinate of $s$.
\end{definition}

\begin{remark}
For any function $\alpha$ defined on the internal edges of $(M, \mathcal T)$, which is called ``angle" and may not satisfy the conditions for defining an angle structure, one can define the combinatorial area $A(d)$ without any difference. By definition, for an internal edge $e$, its combinatorial curvature $k(e)$ is $2\pi$ minus the sum of the angles surrounding it. For the area-curvature function $(A, k)$, one may define \cite{LT} its Euler characteristic function as
$$\chi^{(A,\,k)}(s)=\frac{1}{2\pi}\Big(\sum_{t\in\bigtriangleup}x_{t}(s)A(t)+\sum_{j=1}^m2z_jk(e_j)\Big),$$
where $s$ is expressed in Definition \ref{s-chi(A)}. However, only $k=0$ case is used in our paper.
\end{remark}

Luo-Tillmann (\cite{LT}, Lemma 15) got the following relationship between $\chi^{*}$ and $\chi^{(A)}$:

\begin{lemma}\label{lemma2}
Suppose $(M, \mathcal T)$ admits an angle (or semi-angle) structure $\alpha$ with $(A)$. If $s\in \mathcal C(M,\mathcal T)$, then
\begin{equation}
\chi^{(A)}(s)=\chi^{*}(s)-\frac{1}{2\pi}\sum_{q\in \square}A(q)x_{q}(s),
\end{equation}
where $\square$ is the set of all the normal quadrilaterals of $s$ and $A(q)$ is the combinatorial area of $q$ induced by the angle (or semi-angle) structure $\alpha$.

\end{lemma}


\subsection{Farkas's lemma}\label{sec1}

If $x=(x_{1},...,x_{k})\in \mathbb{R}^{k}$, we use $x>0$ ($x\geq 0$, $x<0$, $x\leq 0$ resp.) to mean that all components of $x_{i}$ are positive (non-negative, negative, non-positive resp.). To approach our main results, we need the following duality result from linear programming, which is known as Farkas's lemma, and can be found, for instance, in \cite{Ziegler}. In the following lemma, vectors in $\mathbb{R}^{k}$ and $\mathbb{R}^{l}$ are considered to be column vectors.

\begin{lemma}\label{lemma}
Let $A$ be a real $k\times l$ matrix, $b\in \mathbb{R}^{k}$, and $\cdot$ be the inner product on $\mathbb{R}^{k}$.
\begin{enumerate}
\item $\{x\in \mathbb{R}^{l}|Ax=b\}\neq \emptyset$ if and only if for all $y\in\mathbb{R}^{k}$ such that $A^{T}y=0$, one has $y\cdot b=0$.
\item $\{x\in \mathbb{R}^{l}|Ax=b,x\geq 0\}\neq \emptyset$ if and only if for all $y\in\mathbb{R}^{k}$ such that $A^{T}y\leq 0$, one has $y\cdot b\leq0$.
\item $\{x\in \mathbb{R}^{l}|Ax=b,x>0\}\neq \emptyset$ if and only if for all $y\in\mathbb{R}^{k}$ such that $A^{T}y\neq0$ and $A^{T}y\leq0$, one has $y\cdot b<0$.
\end{enumerate}
\end{lemma}


\section{The proof of Proposition \ref{main1}}\label{3}
It must be pointed out that in the proof of our main theorems, in addition to ideal ones, we only use truncated type 1-3  tetrahedra and type 1-$k$  polyhedra. Therefore, whenever a \emph{hyperideal polyhedron} (or \emph{partially truncated polyhedron}, \emph{hyperideal tetrahedron}, \emph{partially truncated tetrahedron} resp.) appears later, it always refers to a hyperideal polyhedron of 1-$k$ type (partially truncated polyhedron of 1-$k$ type, hyperideal tetrahedron of 1-3 type, partially truncated tetrahedron of 1-3 type resp.). And from now on, $M$ always refers to a hyperbolic 3-manifold with both cusps and totally geodesic boundaries.

\subsection{The ideal triangulation of $M$}
\label{ideal-triangulation}
\begin{figure}[htbp]
\centering
\includegraphics[scale=0.27]{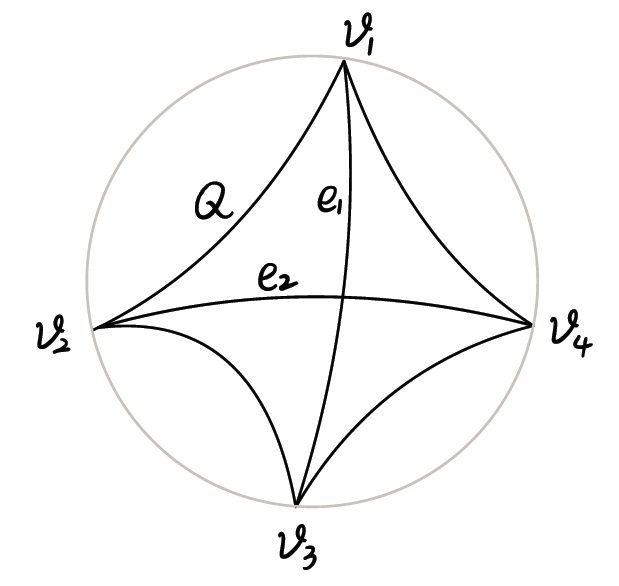}
\caption{a flat ideal tetrahedron}
\label{Fig3}
\end{figure}

\begin{definition}
Let $Q$ be an ideal quadrilateral with four ideal vertices $v_{1}$, $v_{2}$, $v_{3}$, $v_{4}$ in turn. We connect the two pairs of vertices, $v_1$ and $v_3$, $v_2$ and $v_4$, with geodesic segments $e_1$ and $e_2$, respectively. Then the ideal quadrilateral $Q$ with diagonals $e_{1}$ and $e_{2}$ is called a \textbf{flat ideal tetrahedron}, see Figure~\ref{Fig3}.
\end{definition}

By Theorem~\ref{the}, we obtain the following (topological) ideal triangulation of $M$:

\begin{corollary}\label{cor}
Let $M$ be a volume finite, non-compact, complete hyperbolic $3$-manifold with totally
geodesic boundary. Then $M$ admits an ideal triangulation $\mathcal T$ such that each
cell is either an (may be flat) ideal tetrahedron or a partially truncated tetrahedron.
\end{corollary}

\begin{proof}
Using Theorem~\ref{the}, we get a mixed ideal polyhedra decomposition $M$ such that each cell $P$ is either an ideal polyhedron or a partially truncated polyhedron with only one truncated hyperideal vertex.

To obtain the ideal triangulation of $M$, two more steps are required. The first step is to triangulate eacn polyhedron $P$ as follows, and there are two cases to consider:

Case $1$, assume $P$ is a partially truncated polyhedron induced by a convex polyhedron $\tilde{P}$ in $RP^{3}$. Recall a pyramid is a polyhedron with one face (known as the ``base") a $n$-polygon and all the other $n$ faces triangles meeting at a common polygon vertex (known as the ``tip") $v$. Consider the polyhedron $\tilde{P}$, by taking cones at the unique hyperideal vertex $v$, which is called the cone vertex, with each face of $\tilde{P}$ disjoint from $v$, we get a decomposition of $\tilde{P}$ into pyramids. These pyramids have a common tip $v$, and each face of $\tilde{P}$ disjoint from $v$ is the base of a pyramid. For any base polygon $D$ in a pyramid, we arbitrarily pick a vertex $w$ of $D$ and decompose $D$ into triangles by taking cones at $w$. The triangulation of $D$ extends to a triangulation of the pyramid into tetrahedra, and no vertices are added throughout the process. In this way, $\tilde{P}$ is decomposed into a union of convex tetrahedra, each of which has only one common hyperideal vertex $v$ in $RP^{3}$. Then after the truncation at the hyperideal vertex $v$ from the intersection of $\mathbb{H}^3$ with $\tilde{P}$, $P$ has a triangulation whose cells are all partially truncated tetrahedron.

Case $2$, assume $P$ is an ideal polyhedron. By taking cones at any ideal vertex $v$, which is also called cone vertex, we get a decomposition of $P$ into pyramids with tips $v$ and bases the faces of $P$ disjoint from $v$. Similar to case $1$, for any base polygon $D$ in a pyramid, we arbitrarily pick a vertex $w$ of $D$ and decompose $D$ into ideal triangles by taking cones at $w$. The ideal triangulation of $D$ extends to a decomposition of the pyramid into ideal tetrahedra. In this way, $P$ is decomposed into a union of ideal tetrahedra with no vertices are added throughout the process.

Consider a common face $D$ of two adjacent polyhedra $P$ and $P'$. If both $P$ and $P'$ are partially truncated, by the constructions in case 1, the triangulation of $D$ in $P$ matches perfectly with the triangulation of $D$ in $P'$. For the remaining possibilities, at least one of $P$ and $P'$ is an ideal polyhedron, the decompositions of the ideal polygon $D$ from $P$ and $P'$ are not always match together to get a triangulation. Such face $D$ is called a \emph{pillow}.

The second step is to insert flat tetrahedra into these pillows to get a layered triangulation as follows.

Let $v$ and $v'$ be the vertices of a face $D$ where the cone is taken at separately from $P$ and $P'$. Let $v_1,\dots,v_i$ and $\omega_1,\dots,\omega_j$ be the remaining vertices of $D$ which are arranged sequentially along the boundaries of $D$ from $v$ to $v'$, see Figure \ref{Fig6} for example. For each diagonal switch from $v_kv'$ to $v_{k+1}v$, $k=1,..., i-1$, and each diagonal switch from $\omega_kv'$ to $\omega_{k+1}v$, $k=1,...,j-1$, we insert a flat tetrahedron with vertices $\{v, v_i,v_{i+1},v'\}$ and a flat tetrahedron with vertices $\{v,\omega_k,\omega_{k+1},v'\}$ respectively.

 \begin{figure}[htbp]
\centering
\includegraphics[scale=0.55]{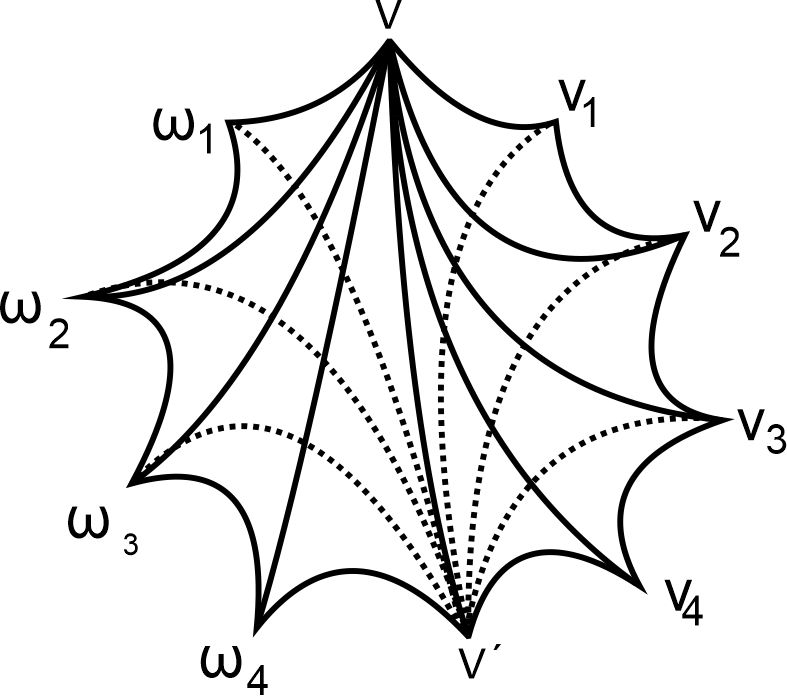}
\caption{triangulating a pillow}
\label{Fig6}
\end{figure}

After the above two steps, we finally get an ideal triangulation $\mathcal T$ of $M$ such that each cell is an (may be flat) ideal tetrahedron or a partially truncated tetrahedron.
\end{proof}

Note at each internal edge in the ideal triangulation $\mathcal T$, the value of its dihedral angle is in $[0,\pi]$, hence we obtain a naturally semi-angle structure on $\mathcal T$. 

\subsection{Sufficient conditions for the existence of angle structures}
\begin{figure}[htbp]
\centering
\includegraphics[scale=0.25]{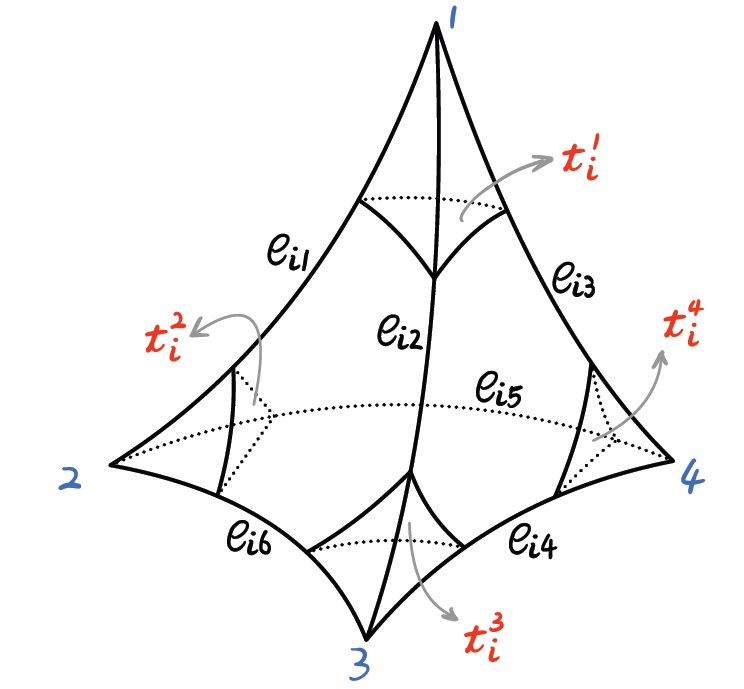}
\caption{the edge labelling and normal triangles in $\sigma_{i}$}
\label{Fig5}
\end{figure}

Let $\mathcal T$ be the ideal triangulation derived in the previous section, and all its (flat) ideal or partially truncated tetrahedra are $\sigma_{1},...,\sigma_n$. Suppose $\alpha$ is an angle structure on $(M, \mathcal T)$. For each $\sigma_i$ and each internal edge $e_{ij}~(1\leq j\leq 6)$ in $\sigma_{i}$, consider the dihedral angle $\alpha_{ij}=\alpha(e_{ij})$ at $e_{ij}$ as a variable of the equations in Farkas' Lemma \ref{lemma}, see Figure \ref{Fig5}.

At the four normal triangles $t_{i}^1$, $...$, $t_{i}^4$ of $\sigma_{i}$, we have the following equations:

\begin{equation}\label{equ4}
  \begin{cases}
\alpha_{i1}+\alpha_{i2}+\alpha_{i3}=a_{i}^{1}\\
\alpha_{i1}+\alpha_{i5}+\alpha_{i6}=a_{i}^{2}\\
\alpha_{i3}+\alpha_{i4}+\alpha_{i5}=a_{i}^{3}\\
\alpha_{i2}+\alpha_{i4}+\alpha_{i6}=a_{i}^{4}.\\
  \end{cases}
\end{equation}

If $\sigma_{i}$ is a partially truncated tetrahedron, then only one of its $a_{i}^{k}$ is less then $\pi$, and the the remaining three are equal to $\pi$. If $\sigma_{i}$ is an (flat) ideal tetrahedron, then all its $a_{i}^{k}$ are equal to $\pi$. The elements of $\{a_{i}^{k}|1\leq i\leq n, 1\leq k\leq 4\}$ are divided into two classes, one of them satisfies $\{a_{i_1}^{k_1}<\pi\}$ with the corresponding normal triangles denoted by $\{t_{i_1}^{k_1}\}$, the other satisfies $\{a_{i_2}^{k_2}=\pi\}$ with the corresponding normal triangles denoted by $\{t_{i_2}^{k_2}\}$. 

Recall $n$ is the number of $3$-simplex in $\mathcal T$ and $m$ is the number of internal edges in $\mathcal T$. At each internal edge $e_j$ of $\mathcal T$, where $1\leq j\leq m$, there holds:

\begin{equation}\label{equ5}
\sum_{i}\alpha_{ij}=2\pi.
\end{equation}
where the sum traverses all the tetrahedra surround the edge $e_j$.

Consider the vector
 \begin{equation}
 (a,b)=(...,a_{i}^{k},...,b_j,...)
 \end{equation}
where $a_{i}^{k}$ ($i=1,...,n$, $k=1,...,4$) comes from (\ref{equ4}) and $b_j=2\pi$ ($j=1,...,m$) comes from (\ref{equ5}). Writing the system of equations (\ref{equ4}) and (\ref{equ5}) in a matrix form as
$$Bx=(a, b)^{T}.$$
Consider the transpose $B^T$ of $B$, which is also the dual of $B$. It has one variable $h_i^k$ for each normal triangle type $t_i^k$ and one variable $z_j$ for each internal edge $e_j$. If we denote the variables corresponding to $B^T $ by $(h, z)=(...,h_{i}^{k},...,...,z_{j},...)$, where $i=1,...,n$, $k=1,...,4$ and $j=1,...,m$, then $B^T(h, z)^T=(...,z_j + h_{i}^{k} + h_{i}^{l},...)^{T}$.

The following useful formula (i.e., Formula (4.9) in \cite{LT}) belongs to Luo-Tillmann:

\begin{equation}\label{equation2}
\begin{aligned}
\frac{1}{\pi}(h, z)\cdot (a, b)&=\chi^{*}(W_{\omega,z})-\chi^{(A)}(W_{\omega,z})\\
&\;\;\;\;+\frac{1}{2\pi}\sum (z_j + h_{i}^{k} + h_{i}^{l})(a_{i}^{k}+a_{i}^{l}-2\pi).
\end{aligned}
\end{equation}
where $W_{\omega,z}=\sum\limits_{i=1}^{n}\omega_i W_{\sigma_{i}}+\sum\limits_{j=1}^{m}z_j W_{e_{j}}$, and the summation runs over all internal edges in the whole triangulation $\mathcal{T}$.

Note the vector $(a,b)$ looks like this
$$(~\underbrace{...~,~a_{i_1}^{k_1},~...~,~\pi,~...~, \pi}_{a},~\underbrace{2\pi,~...~,2\pi}_{b}~),$$
where $a_{i_1}^{k_1}<\pi$. In order to obtain an angle structure on $(M, \mathcal T)$, we only need to prove that there exists a vector $(a,b)=(...,a_{i_1}^{k_1},...,\pi,...,2\pi,...)$ such that
\begin{equation}
\big\{Bx=(a,b)^T, \,x>0\big\}\neq \emptyset.
\end{equation}

Further using the third part of Farkas's lemma, i.e., Lemma \ref{lemma}, we just need to show that for all $(h, z)$ with $B^{T}(h,z)^T\neq 0$ and $B^{T}(h,z)^T\leq 0$, one has $(h,z)\cdot (a,b)< 0$.

\begin{claim}\label{claim}
There exists a vector $(...,a_{i_{1}}^{k_{1}},...,\pi,...,2\pi,...)=(a,b)$ with $a_{i_{1}}^{k_{1}}<\pi$ such that
$$(h,z)\cdot (a,b)< 0$$
for all $(h, z)$ with $B^{T}(h,z)^T\neq 0$ and $B^{T}(h,z)^T\leq 0$.

\end{claim}
\begin{proof}
By Corollary~\ref{cor},  $\mathcal T$ is an ideal triangulation of $M$. Then $(M,\mathcal T)$ admits a semi-angle structure according to the decomposition operation in Section \ref{ideal-triangulation} (when inserting a flat polyhedron, a $0$ angle or $\pi$ angle may occur). Then by the second part of Farkas's Lemma \ref{lemma}, there exists a $\bar{a}_{i_{1}}^{k_{1}}<\pi$ such that $(h,z)\cdot (\bar{a}, b)\leq 0$ for all $B^{T}(h,z)^T\leq 0$, where $(\bar{a}, b)=(...,\bar{a}_{i_{1}}^{k_{1}},...,\pi,...,2\pi,...)$.

According to (\ref{equation2}), as well as the known conditions in Theorem~\ref{main1}, if $a_{i_{1}}^{k_{1}}=\pi$, then

\begin{equation*}\label{equ6}
\begin{aligned}
&\;\;\;\;\frac{1}{\pi}(h, z)\cdot (a, b)\\
&=\chi^{*}(W_{\omega,z})-\chi^{(A)}(W_{\omega,z})+\frac{1}{2\pi}\sum (z_j + h_{i}^{k} + h_{i}^{l})(a_{i}^{k}+a_{i}^{l}-2\pi)\\
&=\chi^{*}(W_{\omega,z})-\chi^{(A)}(W_{\omega,z})\\
&=\chi^{*}(W_{\omega,z})<0.
\end{aligned}
\end{equation*}

Next, we use the continuity method to prove the above, and for this, we introduce a function $F(t)$ defined for all $t\in [0,1]$:
\begin{equation*}
\begin{aligned}
F(t)(h,z)&=(h,z)\cdot (...,\bar{a}_{i_{1}}^{k_{1}}+t(\pi-\bar{a}_{i_{1}}^{k_{1}}),...,\pi,...,2\pi,...)\\
&=\sum(\bar{a}_{i_{1}}^{k_{1}}+t(\pi-\bar{a}_{i_{1}}^{k_{1}}))h_{i_{1}}^{k_{1}}+\sum \pi h_{i_{2}}^{k_{2}}+\sum_{j=1}^{m}2\pi z_{j}.
\end{aligned}
\end{equation*}
Here $(h,z)$ are the variables corresponding to $B^T $, and $h_{i_{2}}^{k_{2}}$ is the variable corresponding to $a_{i}^{k}=\pi$ in the equation (\ref{equ4}).

Thus we have $F(0)\leq 0$, and by Equation (\ref{equ6}), $F(1)< 0$ for all $B^{T}(h,z)^T\neq 0$ and $B^{T}(h,z)^T\leq 0$.

\begin{lemma}
$F(t)(h,z)<0$ for any $t\in (0,1)$, and for any variable $(h,z)$.
\end{lemma}
\begin{proof}
We prove it by contradiction. Assume there is a $t_0\in(0,1)$ and $(h_0,z_0)$ such that
\begin{equation*}
F(t_0)(h_0,z_0)\geq 0.
\end{equation*}
On the one hand $F(0)(h_0,z_0)\leq 0$, derived from the seme-angle structure on ($M, \mathcal T$), and on the other hand, $F(t)(h_0,z_0)$ as a function of $t$ is monotonic because of the linearity of $F(t)(h_0,z_0)$. Then $F(1)(h_0,z_0)\geq 0$, which contradicts the fact that $F(1)(h_0,z_0)<0$. Then we proved the above lemma.

\end{proof}
Finally, if set $a_{i_{1}}^{k_{1}}=\bar{a}_{i_{1}}^{k_{1}}+t_{0}(\pi-\bar{a}_{i_{1}}^{k_{1}})<\pi$, then we have
$$(h,z)\cdot (a,b)=F(t_{0})(h,z)< 0.$$
Hence we get the above Claim \ref{claim}.
\end{proof}

By the third part of Farkas's Lemma \ref{lemma}, $$\big\{Bx=(a,b)^T,x> 0\big\}\neq \emptyset.$$
Hence $(M,\mathcal T)$ admits an angle structure, and Proposition \ref{main1} is proved.

\section{The proof of Theorem \ref{main2}}\label{4}
Let $M$ be a cusped hyperbolic $3$-manifold with totally geodesic boundary and $\mathcal T$ be the ideal triangulation of $M$ as in Corollary~\ref{cor}. Recall that there exists a natural semi-angle structure $\alpha$ on $(M, \mathcal{T})$ by the analysis in Section \ref{3}. For a quadrilateral type normal disk $q$ in $\mathcal T$, the combinatorial area of $q$ by definition is
$$A(q)=\alpha_1+\alpha_2+\alpha_3+\alpha_4-2\pi,$$
where $\alpha_1,..., \alpha_4$ are the dihedral angles assigned by $\alpha$ in the vertices (which are the intersections of $q$ with the four edges of a particular tetrahedron $\sigma_i$) of $q$.

\begin{definition}
A quadrilateral type normal disk $q$ is called \textbf{vertical} if $A(q)=0$.
\end{definition}


From Proposition \ref{main1} and Section \ref{3}, we get the following proposition:
\begin{proposition}\label{prop0}
If there is no such $s \in \mathcal{C}(M, \mathcal{T})$, all of its quadrilateral coordinates are non-negative, all non-vertical quadrilateral coordinates are zero, and at least one quadrilateral coordinate is positive, then $(M, \mathcal{T})$ admits an angle structure.
\end{proposition}

\begin{proof}
Let $\alpha$ be the semi-angle structure on $(M,\mathcal T)$, and $q$ be a quadrilateral type normal disk in a particular tetrahedron $\sigma_i$ of $\mathcal T$.

If $\sigma_i$ is a flat tetrahedron, then $A(q)=0+\pi+0+\pi-2\pi=0$. Hence $q$ is vertical.

If $\sigma_i$ is an ideal or partially truncated tetrahedron, in this case it is not flat, then the combinatorial area
\begin{equation*}
\begin{split}
A(q)&=\alpha_1+\alpha_2+\alpha_3+\alpha_4-2\pi\\
&<\alpha_1+\alpha_2+\alpha_3+\alpha_4+2\alpha_5-2\pi \\
&=(\alpha_1+\alpha_2+\alpha_5-\pi)+(\alpha_3+\alpha_4+\alpha_5-\pi)\\
&\leq0,
\end{split}
\end{equation*}
where $\alpha_5$ is the angle of the edge in $\sigma_i$ which faces $q$, and is positive as a real dihedral angle induced from the hyperbolic structure.

According to the conditions in the proposition, for any  $s\in \mathcal C(M,\mathcal T)$ with all quadrilateral coordinates non-negative and at least one quadrilateral coordinate positive, there
exist quadrilateral type normal disk $q$ such that $q$ is non-vertical and $x_q>0$. Note $\chi^{(0)}(s)=0$, hence by Lemma~\ref{lemma2},
\begin{equation*}
\chi^{*}(s)=\chi^{*}(s)-\chi^{(0)}(s)=\frac{1}{2\pi}\sum_{q\in \square}A(q)x_{q}(s)<0.
\end{equation*}
Then by Theorem~\ref{main1}, we get an angle structure on $\mathcal{T}$.
\end{proof}

Recall $\mathcal{C}(M, \mathcal{T})\subset \mathbb{R}^7$ is the solution space of the compatibility equations (\ref{equation1}). Denote $\mathcal{N}(M, \mathcal{T})$ by a subset of $\mathcal{C}(M, \mathcal{T})$ whose element has all non-negative integer coordinates. With a similar argument of Hodgson-Rubinstein-Segerman~\cite{HRS}, we have:

\begin{proposition}\label{pro1}
If there is no such $x \in \mathcal{N}(M, \mathcal{T})$, all of its quadrilateral coordinates are non-negative, all non-vertical quadrilateral coordinates are zero, and at least one quadrilateral coordinate is positive, then $(M, \mathcal{T})$ admits an angle structure.
\end{proposition}

\begin{proof}
If there is a $x \in \mathcal{C}(M, \mathcal{T})$ satisfying the properties assumed in Proposition \ref{prop0}, then the coordinates of $x$ are all rational numbers since the coefficients in the compatibility equations  (\ref{equation1}) are all integers. So there is a positive integer $N$ such that $Nx\in \mathcal{C}(M, \mathcal{T})$, i.e., $Nx$ has integer coordinates. Because all the quadrilateral coordinates of $x$ are non-negative, so is $Nx$. If all the triangle coordinates of $Nx$ are non-negative, we are done. If some triangle coordinates of $Nx$ are negative, we can always add enough normal copies of peripheral $\partial M$, tori or Klein bottles. Then we can get a $(Nx)'\in \mathcal{C}(M, \mathcal{T})$ whose triangle coordinates are all none-negative. Hence $(Nx)'\in \mathcal{N}(M; \mathcal{T})$, we are done.

\end{proof}

According to Proposition~\ref{pro1}, to get an angle structure on $(M, \mathcal{T})$, we should exclude the cases that ``$x \in \mathcal{N}(M, \mathcal{T})$ with all its quadrilateral coordinates non-negative, all non-vertical quadrilateral coordinates zero, and at least one quadrilateral coordinate positive".

Let $\mathcal{D}$ be the mixed polyhedron decomposition described in Theorem~\ref{the}. Let $\hat{D}$ be the decomposition synthesized by $\mathcal{D}$ with the ideal polyhedra pillows defined in the proof of Corollary~\ref{cor}. Let $\{P_j, j=1,...,k'\}$ and $\{D_{j'},j'=k'+1,...,n'\}$ be the ideal polyhedra (including partially truncated polyhedra) and polygonal pillows respectively in $\hat{D}$, and let $\{\sigma_{i'},i'=1,...,k\}$ and $\{\sigma_{i''},i''=k+1,...,n\}$ be the non-planar 3-simplices and planar tetrahedra respectively in $\mathcal {T}$. Notably, according to Lemma 7.4 of Hodgson-Rubinstein-Segerman~\cite{HRS}, each polygonal pillow represented by $D_{j'}$ can only be a quadrilateral or a hexagon.

\begin{figure}[htbp]
\centering
\includegraphics[scale=0.58]{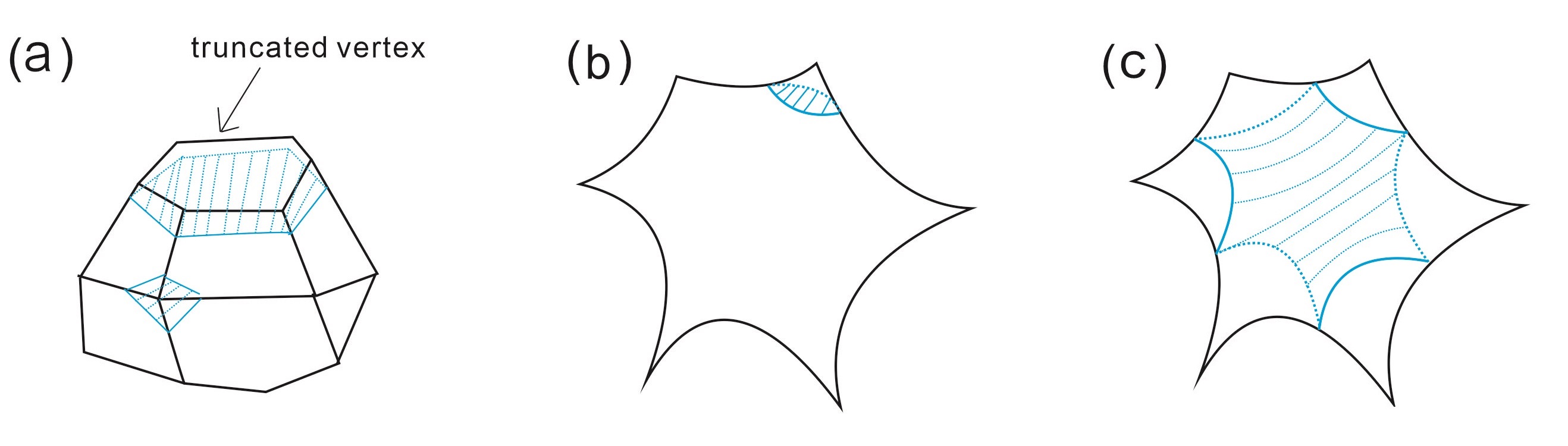}
\caption{$(a)$ is a vertex-linking normal disk; $(b)$ is a vertex-linking bigon in a pillow; $(c)$ is a twisted $6$-gon in a pillow}
\label{Fig7}
\end{figure}

Consider the normal surface $S$ corresponds to a $s \in \mathcal{N}(M, \mathcal{T})$. If $S \cap \sigma_{i'}$ does not contain any normal quadrilaterals for each $i'=1,...,k$, and $S \cap \sigma_{i''}$ contains at least one normal quadrilateral for some $k+1\leq i''\leq n$, then
\begin{enumerate}
\item for each $j=1,...,k'$ so that $S \cap P_j\neq\varnothing$, then $S\cap P_j$ are vertex-linking polygons;
\item for each $j'=k'+1,...,n'$ so that $S \cap D_{j'}\neq\varnothing$, then $S \cap D_{j'}$ are either twisted surfaces or vertex-linking bigons.
\end{enumerate}
See Figure~\ref{Fig7} for example.

The following are the proof of Theorem~\ref{main2}:
\begin{proof}
Let $M$ be a non-compact, volume-finite hyperbolic three-dimensional manifold. Denote $\overline{M}$ by the compact $3$-manifold with boundary, with each torus (or Klein bottle) boundary component corresponds to a cusp of $M$. After subtracting the torus (or Klein bottle) boundary components, $\overline{M}$ is homeomorphic to $M$, and each boundary of $\overline{M}$ is a torus or Klein bottle or a closed surface of high genus. Consider the mixed polyhedral decomposition $\mathcal{D}$ of $M$ described in Theorem~\ref{the}. It determines a graph $\Gamma$ with one vertex for each polyhedron and an edge joining two vertices for each face incident to the two corresponding polyhedra

\begin{lemma}
There exists an ideal triangulation $\mathcal{T}'$ of $M$, which is a subdivision of $\mathcal{D}$, such that $\Gamma$ contains a maximal tree $T$ with the property that the gluing faces corresponding to the edges of $T$ are not polygon pillows in $(M, \mathcal{T}')$.
\end{lemma}
\begin{proof}
For each vertex of degree $1$ in the maximal tree $T$, if the corresponding polyhedron is a truncated one, take the hyperideal vertex as the cone vertex; if the corresponding polyhedron is an ideal one, take the ideal vertex which is not on the face and is represented by an edge of $T$ as the cone vertex. Then remove these vertices and the corresponding edges from $T$, and continue the above operation in the remaining tree until either there is an empty set left, or a single point set left.

If there is an empty set left, then the proof is complete. If there is a single point set left, there are two scenarios to consider:

Case $1$, if the corresponding polyhedron is ideal, then any triangulation of this polyhedron satisfies the condition.

Case $2$, if the corresponding polyhedron is truncated, then take the unique hyperideal vertex of this truncated polyhedron as the cone vertex, and take any triangulation at the other ideal faces. Since each truncated face is shared by two truncated polyhedra, the inherited triangulations of the truncated faces match perfectly. Hence the final triangulation of this truncated polyhedron also satisfies the condition.
\end{proof}

Note that the construction of $\mathcal{T}'$ and $\mathcal{T}$ is consistent as in Corollary~\ref{cor}. Hence $(M; \mathcal{T}')$ also satisfies Proposition \ref{main1}.

For the ideal triangulation $\mathcal{T}'$, denote $\mathcal{D}'$ by the union of $\mathcal{D}$ and polygon pillows generated by $\mathcal{T}'$. Now assume that there exists a $s \in \mathcal{N}(M, \mathcal{T}')$, and without loss of generality, further assume that the corresponding normal surface $S$ (of $s$) is not a Haken sum and contains some twisted disks. By the conditions of Theorem \ref{main2} and the construction of $\mathcal{T}'$, $S$ does not intersect the maximal tree $T$.

Next we show the number of twisted disks in $S$ is even. If there is an odd number of twisted disks of $S$ at some polygon pillow in $\mathcal{D}'$. Denote the edge on $\Gamma$ corresponding to the polygon pillow by $e$. Then $T\cup e$ contains a closed curve $c$ whose intersection with $S$ is odd. Thus, $c$ represents a nontrivial element $[c]$ in the homology group $H_1(\overline{M}; Z_2)$.

Consider the following long exact sequence:
\begin{equation*}
...\rightarrow H_2(\overline{M},\partial \overline{M};Z_2)\rightarrow H_1(\partial \overline{M};Z_2) \stackrel{f}{\rightarrow} H_1(\overline{M};Z_2)\stackrel{g}{\rightarrow} H_1(\overline{M},\partial \overline{M};Z_2)\rightarrow ...
\end{equation*}

Since the closed curve representing any non-trivial element in $H_1(\partial \overline{M};Z_2)$  must have an even geometric intersection number with any closed surface in $M$, $[c]$ is not in the image of the map $f$, and therefore not in the kernel of the map $g$, contradicting the known conditions in Theorem \ref{main2}. Therefore, the number of twisted disks in $S$ is even at all polygon pillows.

If the number of all vertex-linking disks and vertex-linking bigons are even, then $S$ is a double of some surface, contradicting the assumption that $S$ is not a Haken sum. Hence, there must be at least one vertex-linking bigons or vertex-linking disks whose number is odd. Let $S'$ be a copy of all such odd-count vertex-linking bigons or vertex-linking disks.

\begin{claim}
$S'$ is a solution to the compatibility equations.
\end{claim}
\begin{proof}
Consider the normal arc type $\alpha$ at a face $F$ in the triangulation $(M, \mathcal{T}')$.

Case $1$: Both tetrahedra $\sigma_1$, $\sigma_2$ adjacent to $F$ are ideal, or truncated. If $S'$ has a normal disk at $F$, it must be a vertex-linking disk in some ideal polyhedron or partially truncated polyhedron of $\mathcal{D}'$. Since the contributions of $\sigma_1$, $\sigma_2$ to $\alpha$ are the same, it satisfies the compatibility equations.

Case $2$: Both tetrahedra $\sigma_1$, $\sigma_2$ adjacent to $F$ are flat ideal. If $S'$ has a normal disk at $F$, it must be a vertex-linking bigon in some polygonal pillows of $\mathcal{D}'$. Again, the contributions of $\sigma_1$, $\sigma_2$ to $\alpha$ are the same, satisfying the compatibility equations.

Case $3$: One tetrahedron $\sigma_1$ is ideal or truncated, and the other $\sigma_2$ is flat ideal. Since the possible twisted disks on $S$ at $\alpha$ are even, the parity of bigons and linking disks at $\alpha$ is the same. If $S'$ has a normal disk at $F$, it must be a bigon and a linking disk in some polygonal pillows of $\mathcal{D}'$. The contributions of $\sigma_1$, $\sigma_2$ to $\alpha$  are the same, also satisfying the compatibility equations.
\end{proof}

Therefore, $S$ is the Haken sum of $S'$ and some surface which contains at least one twisted disk. This contradicts with the assumption that $S$ is not a Haken sum. Therefore Theorem~\ref{main2} is proved.
\end{proof}




~

\noindent
\noindent

\noindent Huabin Ge, hbge@ruc.edu.cn\\[2pt]
\emph{School of Mathematics, Renmin University of China, Beijing 100872, P. R. China}\\[2pt]

\noindent Longsong Jia, jialongsong@stu.pku.edu.cn\\[2pt]
\emph{School of Mathematical Sciences, Peking University, Beijing, 100871, P.R. China}\\[2pt]

\noindent Faze Zhang, zhangfz201@nenu.edu.cn\\[2pt]
\emph{School of Mathematics and Statistics, Northeast Normal University, Changchun, Jilin, 130024, P.R.China} \\[2pt]

\end{document}